\documentclass[12pt]{article}
\usepackage{amsfonts,epsf,amsmath,amssymb}
\usepackage[numbers]{natbib}
\usepackage{epsfig}
\usepackage{latexsym}
\usepackage{pgf,tikz,subfigure}
\usetikzlibrary{decorations.markings}

\newtheorem{theorem}{\bf Theorem}[section]

\newtheorem{lemma}[theorem]{\bf Lemma}
\newtheorem{proposition}[theorem]{\bf Proposition}

\newcommand{\proof}{\noindent{\bf Proof.\ }}
\newcommand{\qed}{\hfill $\square$ \bigskip}

\begin{document}

\title{PARALLELISM OF STABLE TRACES}

\author{
Jernej Rus \\
Abelium d.o.o.\\
Kajuhova 90, 1000 Ljubljana, Slovenia \\
jernej.rus@gmail.com
}
\date{\today}
\maketitle

\begin{abstract}
A parallel $d$-stable trace is a closed walk which traverses every edge of a graph exactly twice in the same direction and for every vertex $v$, there is no subset $X \subseteq N(v)$ with $1 \leq |N| \leq d$ such that every time the walk enters $v$ from $X$, it also exits to a vertex in $X$. In the past, $d$-stable traces were investigated as a mathematical model for an innovative biotechnological procedure -- self-assembling of polypeptide structures. Among other, it was proven that graphs that admit parallel $d$-stable traces are precisely Eulerian graphs with minimum degree strictly larger than $d$. In the present paper we give an alternative, purely combinatorial proof of this result.
\end{abstract}

\noindent
{\bf Keywords:} Eulerian graph; parallel $d$-stable trace; nanostructure design; self-assembling; polypeptide

\medskip\noindent
{\bf AMS Subject Classification (2010):}  05C45, 
05C85, 
94C15 

\newpage
\section{Introduction}

All graphs considered in this paper will be connected, finite, and simple, that is, without loops and multiple edges. If $v$ is a vertex of a graph $G$, then its degree will be denoted by $d_G(v)$ or $d(v)$ for short if $G$ will be clear from the context. The {\em minimum} and the {\em maximum} degree of $G$ are denoted with $\delta(G)$ and $\Delta(G)$, respectively. A {\em directed graph} is a graph where edges have a direction associated with them. In formal terms a directed graph is a pair $G = (V,A)$, where $V$ is a set of vertices and $A$ is a set of ordered pairs of vertices, called {\em arcs}. A maximal connected subgraph of $G$ is called a {\em component} of $G$, while a vertex which separates two other vertices of the same component is a {\em cutvertex}, and an edge separating its ends is a {\em bridge}. A maximal connected subgraph without a cutvertex is called a {\em block}. Thus, every block of a graph $G$ is either a maximal $2$-connected subgraph, or a bridge (with its ends), or an isolated vertex. A {\em subtree} $T$ of a graph $G$ is a subgraph of $G$ that is also a tree (any pair of vertices $u, v \in V(T) \subseteq V(G)$ are connected by exactly one path in $T$).  For other general terms and concepts from graph theory not recalled here we refer to~\cite{we-1996}. 

A {\em circuit} is a closed walk allowing repetitions of vertices and edges. An {\em Eulerian circuit} in $G$ is a circuit which traverses every edge of $G$ exactly once. $G$ is called {\em Eulerian} if it admits an Eulerian circuit. A {\em double trace} in a graph $G$ is a circuit that traverses every edge exactly twice. For a set of vertices $X \subseteq N(v)$, we say that a double trace $W$ has an {\em $X$-repetition} at vertex $v$ (nontrivial $X$-repetition in~\cite{fi-2013}), if $X$ is nonempty, $X \neq N(v)$, and whenever $W$ comes to $v$ from a vertex in $X$ it also continues to a vertex in $X$. An $X$-repetition (at $v$) is a {\em $d$-repetition} if $|X|=d$ (repetition of {\em order} $d$), see Fig.~\ref{fig:repetition}. Clearly if $W$ has an $X$-repetition at $v$, then it also has an $N(v)\setminus X$-repetition at $v$ ({\em symmetry of repetitions}). We call a double trace without any repetition of order $\leq d$ a {\em $d$-stable trace}. Note, that for every $d' \leq d$, a $d$-stable trace is also a $d'$-stable trace.

\begin{figure}[ht!]
\begin{center}
\begin{tabular}{c c c c}
\subfigure
{
\begin{tikzpicture}[scale=0.5,style=thick]
\fill (2,2) circle (3pt) node[below left]{$u$};
\fill (0,2) circle (2pt);
\fill (4,2) circle (2pt);
\fill (2,0) circle (2pt);
\fill (2,4) circle (3pt);
\draw[decoration={markings, mark=at position 0.25 with {\arrow{>}}},postaction={decorate}] (1.9,4)  .. controls (1.8,1.5) and (2.2,1.5) .. (2.1,4);
\draw[decoration={markings, mark=at position 0.75 with {\arrow{>}}},postaction={decorate}] (1.9,4)  .. controls (1.8,1.5) and (2.2,1.5) .. (2.1,4);
\draw (0,2)--(2,2);
\draw (4,2)--(2,2);
\draw (2,0)--(2,2);
\draw (2.1,3) node[right] {$e$};
\end{tikzpicture}
}
&
\subfigure
{
\begin{tikzpicture}[scale=0.5,style=thick]
\fill (2,2) circle (3pt) node[below left]{$v$};
\fill (0,2) circle (2pt);
\fill (4,2) circle (3pt);
\fill (2,0) circle (2pt);
\fill (2,4) circle (3pt);
\draw[decoration={markings, mark=at position 0.25 with {\arrow{>}}},postaction={decorate}] (1.9,4)  .. controls (1.8,1.7) and (1.7,1.8) .. (4,1.9);
\draw[decoration={markings, mark=at position 0.75 with {\arrow{>}}},postaction={decorate}] (1.9,4)  .. controls (1.8,1.7) and (1.7,1.8) .. (4,1.9);
\draw[decoration={markings, mark=at position 0.25 with {\arrow{>}}},postaction={decorate}] (2.05,4)  .. controls (2,1.7) and (1.7,2) .. (4,2.05);
\draw[decoration={markings, mark=at position 0.75 with {\arrow{>}}},postaction={decorate}] (2.05,4)  .. controls (2,1.7) and (1.7,2) .. (4,2.05);
\draw (0,2)--(2,2);

\draw (2,0)--(2,2);
\end{tikzpicture}
}
&
\subfigure
{
\begin{tikzpicture}[scale=0.5,style=thick]
\fill (2,2) circle (3pt) node[below left]{$v$};
\fill (0,2) circle (2pt);
\fill (4,2) circle (3pt);
\fill (2,0) circle (2pt);
\fill (2,4) circle (3pt);
\draw[decoration={markings, mark=at position 0.25 with {\arrow{<}}},postaction={decorate}] (1.9,4)  .. controls (1.8,1.7) and (1.7,1.8) .. (4,1.9);
\draw[decoration={markings, mark=at position 0.75 with {\arrow{<}}},postaction={decorate}] (1.9,4)  .. controls (1.8,1.7) and (1.7,1.8) .. (4,1.9);
\draw[decoration={markings, mark=at position 0.25 with {\arrow{>}}},postaction={decorate}] (2.05,4)  .. controls (2,1.7) and (1.7,2) .. (4,2.05);
\draw[decoration={markings, mark=at position 0.75 with {\arrow{>}}},postaction={decorate}] (2.05,4)  .. controls (2,1.7) and (1.7,2) .. (4,2.05);
\draw (0,2)--(2,2);
\draw (2,0)--(2,2);
\end{tikzpicture}
}
&
\subfigure
{
\begin{tikzpicture}[scale=0.5,style=thick]
\fill (3,2) circle (3pt);
\fill (2.2,2) node[left]{$w$};
\fill (0,0) circle (3pt);
\fill (3,0) circle (3pt);
\fill (6,0) circle (3pt);
\fill (0,4) circle (3pt);
\fill (3,4) circle (3pt);
\fill (6,4) circle (3pt);
\draw (-0.1,4)  .. controls (2.8,1.45) and (3.2,1.45) .. (6.1,4);
\draw (0.1,4)  .. controls (3.1,1.45) and (2.9,2) .. (2.9,4);
\draw (3.1,4)  .. controls (2.9,1.45) and (3.1,2) .. (5.9,4);
\draw (-0.1,0)  .. controls (2.8,2.55) and (3.2,2.55) .. (6.1,0);
\draw (0.1,0)  .. controls (3.1,2.55) and (2.9,2) .. (2.9,0);
\draw (3.1,0)  .. controls (2.9,2.55) and (3.1,2)  .. (5.9,0);
\end{tikzpicture}
}
\\
\end{tabular}
\caption{Possible $1$-, $2$- and $3$-repetitions at vertices $w$, $u$ and $w$, respectively}
\label{fig:repetition}
\end{center}
\end{figure}

In order to present a mathematical model for the biotechnological procedure from~\cite{gr-2013} graphs that admit $d$-stable traces were characterized in~\cite{fi-2013} (thus generalizing results of Sabidussi~\cite{sa-1977} and Eggleton and Skilton~\cite{eg-1984} about $1$-stable traces and Klav\v zar and Rus~\cite{kl-2013} about $2$-stable traces) as follows:

\begin{proposition}
\label{prop:stable}
{\rm \cite[Proposition~$3.4$]{fi-2013}}
A connected graph $G$ admits a $d$-stable trace if and only if $\delta(G) > d$. 
\end{proposition} 

Let now $W$ be a double trace of a graph $G$. Then every edge $e=uv$ of $G$ is traversed exactly twice. If in both cases $e$ is traversed in the same direction (either both times from $u$ to $v$ or both times from $v$ to $u$) we say that $e$ is a {\em parallel edge} (with respect to $W$). If this is not the case we say that $e$ is an {\em antiparallel edge}. The condition that all the edges of $G$ are of the same type is called a {\em parallelism}. A double trace $W$ is a {\em parallel double trace} if every edge of $G$ is parallel and an {\em antiparallel double trace} if every edge of $G$ is antiparallel. 

By replacing every edge of a graph with two new edges we can quickly prove that every graph (resp.~every Eulerian graph) admits an antiparallel (resp.~parallel) double trace, observation made by several authors, Klav\v zar and Rus in~\cite{kl-2013} among others. While graphs admitting antiparallel $d$-stable traces were thoroughly studied in~\cite{rus-2016}, the characterization of parallel $d$-stable traces was only mentioned as a consequence in~\cite{fi-2013}:

\begin{theorem}
\label{thm:parallel_d}
{\rm \cite[Theorem~$5.4$]{fi-2013}}
A graph $G$ admits a parallel $d$-stable trace if and only if $G$ is Eulerian and $\delta(G) > d$.
\end{theorem}

In the present paper we will give an alternative proof of this result. Instead of graph embeddings (heavily used in~\cite{fi-2013}), our approach to the problem will be purely combinatorial. Characterizing graphs that admit parallel $d$-stable traces also represents a new problem related with forbidden transitions in Eulerian tours of Eulerian graphs (further related problems can be found, among others in~\cite{fl-1990, fl-1991}). 

We can note right away that parallel double traces do not contain $1$-repetitions. Note also that none of the operations that we will use on double traces (concatenations, contractions, deletions, inductive constructions, and reordering) will change the orientation of the edges.

\subsection{Biotechological background}
\label{sec:biotechology}

In $2013$ Gradi\v sar et al.~\cite{gr-2013} presented a novel self-assembly strategy for polypeptide nanostructure design. Their strategy relied on routing a single polypeptide chain consisting of $12$ segments through $6$ edges of the tetrahedron in such a way that every edge was traversed exactly twice. The required mathematical support for the particular case of the tetrahedron and the general case of a polyhedron was already given in~\cite{fi-2013, gr-2013, kl-2013}, where the authors explained that polyhedron $P$ that is composed from a single polymer chain can be naturally represented by a graph $G(P)$ of the polyhedron. Circuits that traverse every edge of $G(P)$ precisely twice, called double traces of $G(P)$, play a key role in modeling the construction process.

The stability of the constructed polyhedra depends on an additional property whether in the double trace the neighborhoods of vertices can be split. The reader interested in the biotehological procedures that motivated our research may also consult the references~\cite{gr-2011, ko-2015}, where the authors also exposed the use of parallel $d$-stable traces.

\section{Graphs admitting parallel $2$-stable traces}
\label{sec:parallel}

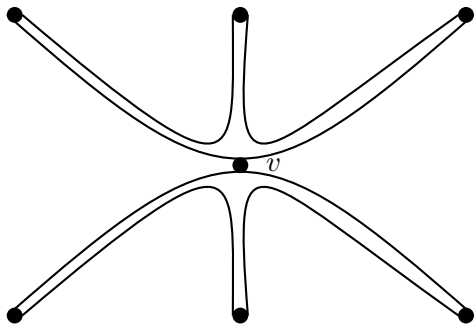
The first mathematical model for the biotechnological procedure from~\cite{gr-2013}, introduced in~\cite{kl-2013}, stated that a polyhedral graph $P$ can be realized by interlocking pairs of polypeptide chains if its corresponding graph $G(P)$ contains a $2$-stable trace. Two important deficiencies of this model were later found in~\cite{fi-2013}: $(i)$ it does not account for vertices of degree $\le 2$, and $(ii)$ it does not successfully model vertices of degree $\ge 6$ (because a polyhedron could split into two parts in a vertex of degree $\ge 6$, as can be seen at Fig.~\ref{fig:repetition} and therefore the structure would not be stable). Since until now, a construction of a polyhedron whose graph would have such properties, has not yet been attempted, we first study parallel $2$-stable traces in this section.

To make the arguments in this section more transparent, we explain how the reader can graphically imagine $1$-repetitions and $2$-repetitions in double traces. We say that a double trace contains a $1$-repetition if it has an immediate succession of an edge $e$ by its antiparallel copy. If $v$ is a vertex of a graph $G$ with a double trace $W$ and $u$ and $w$ are two different neighbors of $v$, then we can say that $W$ contains a $2$-repetition (through) $v$ if the vertex sequence $u\rightarrow v\rightarrow w$ appears twice in $W$ in any direction ($u\rightarrow v\rightarrow w$ or $w\rightarrow v\rightarrow u$), see Fig.~\ref{fig:repetition}.

We will need the next lemma in the proof of Theorem~\ref{thm:parallel}.

\begin{lemma}
\label{lem:contraction}
Let $G$ be a graph and let $T$ a subtree of $G$ such that every vertex $v  \in V(G) \setminus V(T)$ has at most one neighbor in $T$. Construct a graph $G'$ from $G$ by contracting $T$ into a single vertex $t$. If $G$ admits a $2$-stable trace $W$ then $G'$ admits a $2$-stable trace $W'$ that traverses edges from $E(G) \cap E(G')$ in the same direction as $W$.
\end{lemma}

\proof
Suppose that the graph $G$ admits a $2$-stable trace $W$. Construct a double trace $W'$ from $W$ as follows. Start in an arbitrary vertex of $V(G) \setminus V(T)$ and follow $W$. Let $a=xy$ be an arc of $W$ that we are currently traversing on our walk along $W$. If $x,y \in V(G) \setminus V(T)$, then we put $xy$ into $W'$ so that the order of arcs from $W$ is preserved. If $x \in V(T)$ and $y \notin V(T)$ then we put $ty$ in $W'$ instead of $a$. Similarly, we replace arcs where $x \notin V(T)$ and $y \in V(T)$  with $xt$. Finally, the occurrences of the arcs from $T$ are ignored in $W'$. 

We claim that $W'$ is a $2$-stable trace of $G'$. Since every edge is traversed twice in $W$, every edge is traversed twice in $W'$. Hence $W'$ is a double trace. If $W'$ is not a $2$-stable trace, there exists a vertex $x \in V(G')$ such that $W'$ has a $1$-repetition or a $2$-repetition at $x$. Denote the neighborhood of vertex $t$ in $G'$ with $N(t)$. We have to consider three cases.

\vspace{5pt}
\noindent
\textbf{Case 1:} $x \notin N(t)$.

\vspace{5pt}
\noindent
It is clear from the construction that if $W'$ had a $1$-repetition or a $2$-repetition at $x$, then $W$ would have a $1$-repetition or a $2$-repetition at $x$, a contradiction. 

\vspace{5pt}
\noindent
\textbf{Case 2:} $x \in N(t)$.

\vspace{5pt}
\noindent
It is again clear from the construction that if $W'$ had a $1$-repetition $yxy$ or a $2$-repetition $yxz$, where $y,z \neq t$, then $W$ would have a $1$-repetition or a $2$-repetition at $x$, a contradiction.

Assume first that $W'$ has a $1$-repetition $txt$. It follows that $W$ should contain $hxg$, where $h, g \in T$. Since every vertex in $V(G) \setminus V(T)$ has at most one neighbor in $T$, $h = g$. Therefore $W$ should contain a $1$-repetition $hxh$, a contradiction.

Assume next that $W'$ has a $2$-repetition $txy$ for some neighbor $y$ of $x$. It follows that $W$ should contain $hxy$ and $gxy$, where $h, g \in T$. Since every vertex in $V(G) \setminus V(T)$ has at most one neighbor in $T$, $h = g$. Therefore $W$ should contain a $2$-repetition $hxy$, a contradiction.

\vspace{5pt}
\noindent
\textbf{Case 3:} $x = t$.

\vspace{5pt}
\noindent
Assume first that $W'$ has a $1$-repetition $yty$ for some neighbor $y$ of $t$. It follows that $W$ should contain $yhAhy$, where $h$ is the unique neighbor of $y$ in $T$ and $A$ is a circuit in $T$. Since $T$ is a tree, the only possibility that circuit appears in a part of a double trace $W$ that is completely included in $T$ is with a $1$-repetition, a contradiction.

Assume next that $W'$ has a $2$-repetition $ytz$ for some neighbors $y$ and $z$ of $t$. It follows that $W$ should contain $yhBgz$ and $yhCgz$, where $h$ is a unique neighbor of $y$ in $T$, $g$ is a unique neighbor of $z$ in $T$, while $B$ and $C$ are $hg$-paths in $T$. Considering the fact that in a tree any two vertices are connected with a unique path, we can argue that $B = C$ and therefore that then $W$ should have a $2$-repetition ($1$-repetition if $h = g$), a contradiction.

\vspace{5pt}

We have thus proved that $W'$ is a $2$-stable trace in $G'$. During the construction of $W'$ we did not change the direction of any arc from $W$. 
\qed

Note that Lemma~\ref{lem:contraction} is, by repetition of the procedure described above, also true for forests (any number of disjoint subtrees).

The following was proven in~\cite{kl-2013}, where it was also observed that a graph $G$ admits a parallel double trace if and only if $G$ is Eulerian.

\begin{proposition}
\label{prp:parallel}
{\rm \cite[Proposition~$5.4$]{kl-2013}}
A connected graph $G$ admits a parallel $1$-stable trace if and only if $G$ is Eulerian. 
\end{proposition}

\proof
Parallelism of any stable trace of a graph $G$ implies that all the vertices of $G$ are of even degree and traversing an arbitrary Eulerian circuit of $G$ twice in the same direction constructs a parallel $1$-stable trace.
\qed

We next prove Theorem~\ref{thm:parallel} about parallel $2$-stable traces and then use it in Section~\ref{sec:parallel_d} to present an alternative proof of Theorem~\ref{thm:parallel_d}.

\begin{theorem}
\label{thm:parallel}
A graph $G$ admits a parallel $2$-stable trace if and only if $G$ is Eulerian and $\delta(G) > 2$.
\end{theorem}

Note that for Eulerian graphs the constraint on the minimal degree of a graph from Theorem~\ref{thm:parallel} is equivalent to $\delta(G) \geq 4$.

\vspace{5pt}

\proof
Suppose that a graph $G$ admits a parallel $2$-stable trace. By definition, every $2$-stable trace is a $1$-stable trace. Thus by Proposition~\ref{prp:parallel}, $G$ is Eulerian and hence by Proposition~\ref{prop:stable} we infer that $\delta(G) \geq 4$. 

For the converse assume that $G$ fulfills the conditions of the theorem. We proceed by induction on $\Delta = \Delta(G)$. 

Let $\Delta = 4$. Then $\delta(G) = \Delta(G) = 4$. By Proposition~\ref{prp:parallel}, $G$ admits a parallel $1$-stable trace $W'$. If $W'$ is not already a $2$-stable trace, $W'$ contains at least one $2$-repetition. We proceed with the second induction on the number $k$ of vertices where $W'$ has $2$-repetitions. Let $k \geq 1$ and let $v$ be one of the vertices where $W'$ has a $2$-repetition. If a $1$-stable trace $W'$ has a $2$-repetition through $v$, where $v$ is a vertex with $d_G(v) = 4$, then it is not difficult to see that $W'$ has two $2$-repetitions through $v$. Let $v_1$, $v_2$, $v_3$, and $v_4$ be the neighbors of $v$. Without loss of generality, we can assume that $A = v_1 \rightarrow v \rightarrow v_2$ is the first and $B = v_3 \rightarrow v \rightarrow v_4$ is the second $2$-repetition through $v$ in $W'$. That means that sequences $A$ and $B$ appear twice in $W'$. Because $W'$ is a parallel $1$-stable trace, there are only two possibilities how occurrences of $A$ and $B$ are arranged in $W'$. These possibilities are $AABB$ (Fig.~\ref{fig:construction1}, left) and $ABAB$ (Fig.~\ref{fig:construction2}, left). Note that we left out all the other vertices in Figs.~\ref{fig:construction1} and~\ref{fig:construction2}.

In the first case we construct a double trace $W$ from $W'$ in $G$ as follows. We start in an arbitrary vertex of $V(G) \setminus \{v\}$ and follow $W'$. Let $a=xy$ be an arc of $W'$ that we are currently traversing on our walk along $W'$. If $x,y \in V(G) \setminus \{v,v_1,v_2,v_3,v_4\}$, then we put $xy$ into $W$ so that the order of arcs from $W'$ is preserved. Put one occurrence of $v_1 \rightarrow v \rightarrow v_2$ and one occurrence of $v_3 \rightarrow v \rightarrow v_4$ in $W$ as well. Replace the remaining occurrence of $v_1 \rightarrow v \rightarrow v_2$ with $v_1 \rightarrow v \rightarrow v_4$ and the remaining occurrence of $v_3 \rightarrow v \rightarrow v_4$ with $v_3 \rightarrow v \rightarrow v_2$, so that $W$ stays connected, see Fig.~\ref{fig:construction1}, right. 

\begin{figure}[ht!]
\begin{center}
\begin{tikzpicture}[scale=0.8,style=thick]
\fill (2,2) circle (3pt) node[above right]{$v$};
\fill (0,2) circle (3pt) node[above]{$v_4$};
\fill (4,2) circle (3pt) node[right]{$v_2$};
\fill (2,0) circle (3pt) node[right]{$v_3$};
\fill (2,4) circle (3pt) node[right]{$v_1$};
\draw[decoration={markings, mark=at position 0.25 with {\arrow{>}}},postaction={decorate}] (1.9,4)  .. controls (1.8,1.7) and (1.7,1.8) .. (4,1.9);
\draw[decoration={markings, mark=at position 0.75 with {\arrow{>}}},postaction={decorate}] (1.9,4)  .. controls (1.8,1.7) and (1.7,1.8) .. (4,1.9);
\draw[decoration={markings, mark=at position 0.25 with {\arrow{>}}},postaction={decorate}] (2.05,4)  .. controls (2,1.7) and (1.7,2) .. (4,2.05);
\draw[decoration={markings, mark=at position 0.75 with {\arrow{>}}},postaction={decorate}] (2.05,4)  .. controls (2,1.7) and (1.7,2) .. (4,2.05);
\draw[decoration={markings, mark=at position 0.25 with {\arrow{>}}},postaction={decorate}] (1.9,0)  .. controls (2,2.3) and (2.3,2) .. (0,1.9);
\draw[decoration={markings, mark=at position 0.75 with {\arrow{>}}},postaction={decorate}] (1.9,0)  .. controls (2,2.3) and (2.3,2) .. (0,1.9);
\draw[decoration={markings, mark=at position 0.25 with {\arrow{>}}},postaction={decorate}] (2.05,0)  .. controls (2.2,2.3) and (2.3,2.2) .. (0,2.05);
\draw[decoration={markings, mark=at position 0.75 with {\arrow{>}}},postaction={decorate}] (2.05,0)  .. controls (2.2,2.3) and (2.3,2.2) .. (0,2.05);
\draw[decoration={markings, mark=at position 0.5 with {\arrow{>}}},postaction={decorate}] (4,2.05)  .. controls (5,3) and (2.5,6) .. (1.9,4);
\draw[decoration={markings, mark=at position 0.5 with {\arrow{>}}},postaction={decorate}] (4,1.9)  .. controls (5,1) and (2.5,-2) .. (2.05,0);
\draw[decoration={markings, mark=at position 0.5 with {\arrow{>}}},postaction={decorate}] (0,2.05)  .. controls (-2,1.45) and (1.5,-2) .. (1.9,0);
\draw[decoration={markings, mark=at position 0.5 with {\arrow{>}}},postaction={decorate}] (0,1.9)  .. controls (-2,2.5) and (1.5,6) .. (2.05,4);

\fill (6,2) circle(0pt) node{$\Longrightarrow$};

\fill (10,2) circle (4pt) node[above right]{$v$};
\fill (8,2) circle (3pt) node[above]{$v_4$};
\fill (12,2) circle (3pt) node[right]{$v_2$};
\fill (10,0) circle (3pt) node[right]{$v_3$};
\fill (10,4) circle (3pt) node[right]{$v_1$};
\draw[decoration={markings, mark=at position 0.25 with {\arrow{>}}},postaction={decorate}] (9.9,4)  .. controls (9.9,1.7) and (10.2,2) .. (8,2.05);
\draw[decoration={markings, mark=at position 0.72 with {\arrow{>}}},postaction={decorate}] (9.9,4)  .. controls (9.9,1.7) and (10.2,2) .. (8,2.05);
\draw[decoration={markings, mark=at position 0.28 with {\arrow{>}}},postaction={decorate}] (10.05,4)  .. controls (10,1.7) and (9.7,2) .. (12,2.05);
\draw[decoration={markings, mark=at position 0.78 with {\arrow{>}}},postaction={decorate}] (10.05,4)  .. controls (10,1.7) and (9.7,2) .. (12,2.05);
\draw[decoration={markings, mark=at position 0.25 with {\arrow{>}}},postaction={decorate}] (9.9,0)  .. controls (9.8,2) and (9.7,2) .. (12,1.9);
\draw[decoration={markings, mark=at position 0.75 with {\arrow{>}}},postaction={decorate}] (9.9,0)  .. controls (9.8,2) and (9.7,2) .. (12,1.9);
\draw[decoration={markings, mark=at position 0.25 with {\arrow{>}}},postaction={decorate}] (10.05,0)  .. controls (10,2.1) and (10.3,2) .. (8,1.9);
\draw[decoration={markings, mark=at position 0.75 with {\arrow{>}}},postaction={decorate}] (10.05,0)  .. controls (10,2.1) and (10.3,2) .. (8,1.9);
\draw[decoration={markings, mark=at position 0.5 with {\arrow{>}}},postaction={decorate}] (12,2.05)  .. controls (13,3) and (10.5,6) .. (9.9,4);
\draw[decoration={markings, mark=at position 0.5 with {\arrow{>}}},postaction={decorate}] (12,1.9)  .. controls (13,1) and (10.5,-2) .. (10.05,0);
\draw[decoration={markings, mark=at position 0.5 with {\arrow{>}}},postaction={decorate}] (8,2.05)  .. controls (6,1.45) and (9.5,-2) .. (9.9,0);
\draw[decoration={markings, mark=at position 0.5 with {\arrow{>}}},postaction={decorate}] (8,1.9)  .. controls (6,2.5) and (9.5,6) .. (10.05,4);
\end{tikzpicture}
\end{center}
\caption{Removing $2$-repetition through $v$ (case $AABB$)}
\label{fig:construction1}
\end{figure}

We construct $W$ analogously in the second case, see Fig.~\ref{fig:construction2}, right.

\begin{figure}[ht!]
\begin{center}
\begin{tikzpicture}[scale=0.8,style=thick]
\fill (2,2) circle (3pt) node[above right]{$v$};
\fill (0,2) circle (3pt) node[above right]{$v_4$};
\fill (4,2) circle (3pt) node[above]{$v_2$};
\fill (2,0) circle (3pt) node[right]{$v_3$};
\fill (2,4) circle (3pt) node[right]{$v_1$};
\draw[decoration={markings, mark=at position 0.25 with {\arrow{>}}},postaction={decorate}] (1.9,4)  .. controls (1.8,1.7) and (1.7,1.8) .. (4,1.9);
\draw[decoration={markings, mark=at position 0.75 with {\arrow{>}}},postaction={decorate}] (1.9,4)  .. controls (1.8,1.7) and (1.7,1.8) .. (4,1.9);
\draw[decoration={markings, mark=at position 0.25 with {\arrow{>}}},postaction={decorate}] (2.05,4)  .. controls (2,1.7) and (1.7,2) .. (4,2.05);
\draw[decoration={markings, mark=at position 0.75 with {\arrow{>}}},postaction={decorate}] (2.05,4)  .. controls (2,1.7) and (1.7,2) .. (4,2.05);
\draw[decoration={markings, mark=at position 0.25 with {\arrow{>}}},postaction={decorate}] (1.9,0)  .. controls (2,2.3) and (2.3,2) .. (0,1.9);
\draw[decoration={markings, mark=at position 0.75 with {\arrow{>}}},postaction={decorate}] (1.9,0)  .. controls (2,2.3) and (2.3,2) .. (0,1.9);
\draw[decoration={markings, mark=at position 0.25 with {\arrow{>}}},postaction={decorate}] (2.05,0)  .. controls (2.2,2.3) and (2.3,2.2) .. (0,2.05);
\draw[decoration={markings, mark=at position 0.75 with {\arrow{>}}},postaction={decorate}] (2.05,0)  .. controls (2.2,2.3) and (2.3,2.2) .. (0,2.05);
\draw[decoration={markings, mark=at position 0.15 with {\arrow{>}}},postaction={decorate}] (4,2.05)  .. controls (5.5,2.4) and (2.5,-2) .. (2.05,0);
\draw[decoration={markings, mark=at position 0.5 with {\arrow{>}}},postaction={decorate}] (0,2.05)  .. controls (-1,3) and (1.5,6) .. (1.9,4);
\draw[decoration={markings, mark=at position 0.5 with {\arrow{>}}},postaction={decorate}] (4,1.9)  .. controls (4.7,0.8) and (2.5,-2) .. (1.9,0);
\draw[decoration={markings, mark=at position 0.5 with {\arrow{>}}},postaction={decorate}] (0,1.9)  .. controls (-1.8,2.5) and (2.2,7) .. (2.05,4);

\fill (6,2) circle(0pt) node{$\Longrightarrow$};

\fill (10,2) circle (4pt) node[above right]{$v$};
\fill (8,2) circle (3pt) node[above right]{$v_4$};
\fill (12,2) circle (3pt) node[above]{$v_2$};
\fill (10,0) circle (3pt) node[right]{$v_3$};
\fill (10,4) circle (3pt) node[right]{$v_1$};
\draw[decoration={markings, mark=at position 0.25 with {\arrow{>}}},postaction={decorate}] (9.9,4)  .. controls (9.9,1.7) and (10.2,1.8) .. (8,1.9);
\draw[decoration={markings, mark=at position 0.75 with {\arrow{>}}},postaction={decorate}] (9.9,4)  .. controls (9.9,1.7) and (10.2,1.8) .. (8,1.9);
\draw[decoration={markings, mark=at position 0.28 with {\arrow{>}}},postaction={decorate}] (10.05,4)  .. controls (10.1,1.7) and (9.8,1.8) .. (12,1.9);
\draw[decoration={markings, mark=at position 0.78 with {\arrow{>}}},postaction={decorate}] (10.05,4)  .. controls (10.1,1.7) and (9.8,1.8) .. (12,1.9);
\draw[decoration={markings, mark=at position 0.25 with {\arrow{>}}},postaction={decorate}] (9.9,0)  .. controls (9.8,2) and (9.7,2) .. (12,2.05);
\draw[decoration={markings, mark=at position 0.75 with {\arrow{>}}},postaction={decorate}] (9.9,0)  .. controls (9.8,2) and (9.7,2) .. (12,2.05);
\draw[decoration={markings, mark=at position 0.25 with {\arrow{>}}},postaction={decorate}] (10.05,0)  .. controls (10,2.1) and (10.3,2.2) .. (8,2.05);
\draw[decoration={markings, mark=at position 0.75 with {\arrow{>}}},postaction={decorate}] (10.05,0)  .. controls (10,2.1) and (10.3,2.2) .. (8,2.05);
\draw[decoration={markings, mark=at position 0.15 with {\arrow{>}}},postaction={decorate}] (12,2.05)  .. controls (13.5,2.4) and (10.5,-2) .. (10.05,0);
\draw[decoration={markings, mark=at position 0.5 with {\arrow{>}}},postaction={decorate}] (8,2.05)  .. controls (7,3) and (9.5,6) .. (9.9,4);
\draw[decoration={markings, mark=at position 0.5 with {\arrow{>}}},postaction={decorate}] (12,1.9)  .. controls (12.7,0.8) and (10.5,-2) .. (9.9,0);
\draw[decoration={markings, mark=at position 0.5 with {\arrow{>}}},postaction={decorate}] (8,1.9)  .. controls (6.2,2.5) and (10.2,7) .. (10.05,4);
\end{tikzpicture}
\end{center}
\caption{Removing $2$-repetition through $v$ (case $ABAB$)}
\label{fig:construction2}
\end{figure}

We claim that in both cases $W$ is a parallel $1$-stable trace of $G$ with at least one vertex with $2$-repetition less than $W'$. Note first that any edge $e = xy$ that appears in $W$ (arcs $xy$ or $yx$ appears in $W$) has its unique corresponding edge $e'$ in $W'$. Any edge $e = xy$ in $W$, where $x \neq v$ and $y \neq v$, is traversed twice in the same direction in $W$ because it is traversed twice in the same direction in $W'$. Four remaining edges ($vv_1$, $vv_2$, $vv_3$, and $vv_4$) are traversed twice in the same direction by construction. Hence $W$ is a parallel double trace. It is also clear from the construction that $W$ is a $1$-stable trace. Finally we need to verify that $W$ has at least one vertex with $2$-repetition less than $W'$. Let $x$ be an arbitrary vertex of $G$ in which $W$ has a $2$-repetition. We have to consider three cases.

\vspace{5pt}
\noindent
\textbf{Case 1:} $x \notin \{v,v_1,v_2,v_3,v_4\}$.

\vspace{5pt}
\noindent
It is clear from the construction that if $W$ has a $2$-repetition through $x$, then also $W'$ has a $2$-repetition through $x$. 

\vspace{5pt}
\noindent
\textbf{Case 2:} $x \in \{v_1,v_2,v_3,v_4\}$.

\vspace{5pt}
\noindent
It is again clear from the construction that if $W$ has a $2$-repetition $yxz$, where $y,z \neq v$, then also $W'$ has a $2$-repetition through $x$.

Similarly, if $W$ has a $2$-repetition $vxy$ for some neighbor $y$ of $x$, then also $W'$ has a $2$-repetition through $x$ since the order of arcs adjacent to $\{v_1,v_2,v_3,v_4\}$ did not change in $W$.

\vspace{5pt}
\noindent
\textbf{Case 3:} $x = v$.

\vspace{5pt}
\noindent
The $1$-stable trace $W'$ had a $2$-repetition (two $2$-repetitions to be more accurate) through $v$ but during the construction of $W$ we manage to remove them both.

\vspace{5pt}

We have thus constructed a $1$-stable trace $W$ which have at least one vertex with $2$-repetition less than $W'$. Hence, it follows by induction assumption that any $4$-regular graph admits a parallel $2$-stable trace.

Assume now that $\Delta \geq 6$ and that any graph $H$ with $\Delta(H) < \Delta$ that fulfills the conditions of Theorem~\ref{thm:parallel} admits a parallel $2$-stable trace. We have to again consider two cases.

\vspace{5pt}
\noindent
\textbf{Case 1:} $\Delta \equiv 2 \pmod{4}$.

\vspace{5pt}
\noindent
Construct the graph $G'$ from $G$ as follows. For every vertex $v$ of degree $\Delta$ (temporary denote its neighbors with $v_1, \ldots , v_{\Delta}$) repeat the following procedure. Remove $v$ from $G$. Add two new vertices $v'$ and $v''$, connect them by an edge, connect $v'$ with $v_1, \ldots , v_{\frac{\Delta}{2}}$, and connect $v''$ with the remaining neighbors of $v$, see Fig.~\ref{fig:construction3}.

\begin{figure}[ht!]
\begin{center}
\subfigure[$G$]
{
\begin{tikzpicture}[scale=1,style=thick]
\fill (1.5,4) circle (2pt) node[above]{$v$};
\fill (0,2) circle (2pt) node[below left]{$v_1$};
\draw (1.5,2) node {$\ldots$};
\fill (3,2) circle (2pt) node[below right]{$v_{\Delta}$};
\fill (-0.2,0) circle (2pt);
\fill (0.2,0) circle (2pt);
\fill (2.8,0) circle (2pt);
\fill (3,0) circle (2pt);
\fill (3.2,0) circle (2pt);
\draw (1.5,4)--(0,2);
\draw (1.5,4)--(3,2);
\draw (0,2)--(-0.2,0);
\draw (0,2)--(0.2,0);
\draw (3,2)--(2.8,0);
\draw (3,2)--(3,0);
\draw (3,2)--(3.2,0);
\end{tikzpicture}
}
\subfigure[$G'$]
{
\begin{tikzpicture}[scale=1,style=thick]
\draw [dashed] (1.25,4)--(2,4);
\fill (1.25,4) circle (2pt) node[above]{$v'$};
\fill (2,4) circle (2pt) node[above]{$v''$};
\fill (-0.5,2) circle (2pt) node[below left]{$v_1$};
\draw (0.45,2) node {$\ldots$};
\fill (1.4,2) circle (2pt) node[below]{$v_{\frac{\Delta}{2}}$};
\fill (2.6,2) circle (2pt) node[below]{$v_{\frac{\Delta}{2}+1}$};
\draw (3.05,2) node {$\ldots$};
\fill (3.5,2) circle (2pt) node[below right]{$v_{\Delta}$};
\fill (-0.7,0) circle (2pt);
\fill (-0.3,0) circle (2pt);
\fill (3.3,0) circle (2pt);
\fill (3.5,0) circle (2pt);
\fill (3.7,0) circle (2pt);
\draw (1.25,4)--(-0.5,2);
\draw (1.25,4)--(1.4,2);
\draw (2,4)--(2.6,2);
\draw (2,4)--(3.5,2);
\draw (-0.5,2)--(-0.7,0);
\draw (-0.5,2)--(-0.3,0);
\draw (3.5,2)--(3.3,0);
\draw (3.5,2)--(3.5,0);
\draw (3.5,2)--(3.7,0);
\end{tikzpicture}
}
\end{center}
\caption{Construction from the proof of Theorem~\ref{thm:parallel} for the case $\Delta \equiv 2 \pmod{4}$}
\label{fig:construction3}
\end{figure}

Note that in $G'$ all except the newly added vertices are of the same degree as in $G$, while $d_{G'}(v') = \frac{\Delta}{2} + 1$ and $d_{G'}(v'') = \frac{\Delta}{2} + 1$ (the last two statements are true for all new vertices). It follows that $\Delta(G') < \Delta$. Since $\Delta \ge 6$, we then also infer that $\delta(G')\ge 4$. Because $\Delta \equiv 2 \pmod{4}$, the degrees $d_{G'}(v') = d_{G'}(v'') = \frac{\Delta}{2} + 1$ are even, hence $G$ is Eulerian and by the induction assumption on $\Delta$, the graph $G'$ admits a parallel $2$-stable trace. If we use a path containing vertices $v'$ and $v''$ as subtree $T$, it follows from a repeated application of Lemma~\ref{lem:contraction} that $G$ admits a parallel $2$-stable trace.

\vspace{5pt}
\noindent
\textbf{Case 2:} $\Delta \equiv 0 \pmod{4}$.

\vspace{5pt}
\noindent
Construct the graph $G'$ from $G$ as follows. For every vertex $v$ of degree $\Delta$ (temporary denote its neighbors with $v_1, \ldots , v_{\Delta}$) repeat the following procedure. Remove $v$ from $G$, and add three new vertices $v'$, $v''$, and $v'''$. Connect $v''$ with $v'$ and $v'''$ by an edge, connect $v'$ with $v_1, \ldots , v_{\frac{\Delta}{2} - 1}$, connect $v''$ with $v_{\frac{\Delta}{2}}$ and $v_{\frac{\Delta}{2} + 1}$,  and connect $v'''$ with the remaining neighbors of $v$, see Fig.~\ref{fig:construction4}. 

\begin{figure}[ht!]
\begin{center}
\subfigure[$G$]
{
\begin{tikzpicture}[scale=1,style=thick]
\fill (1.5,4) circle (2pt) node[above]{$v$};
\fill (0,2) circle (2pt) node[below left]{$v_1$};
\draw (1.5,2) node {$\ldots$};
\fill (3,2) circle (2pt) node[below right]{$v_{\Delta}$};
\fill (-0.2,0) circle (2pt);
\fill (0.2,0) circle (2pt);
\fill (2.8,0) circle (2pt);
\fill (3,0) circle (2pt);
\fill (3.2,0) circle (2pt);
\draw (1.5,4)--(0,2);
\draw (1.5,4)--(3,2);
\draw (0,2)--(-0.2,0);
\draw (0,2)--(0.2,0);
\draw (3,2)--(2.8,0);
\draw (3,2)--(3,0);
\draw (3,2)--(3.2,0);
\end{tikzpicture}
}
\subfigure[$G'$]
{
\begin{tikzpicture}[scale=1,style=thick]
\draw [dashed] (1,4)--(2.5,4.5);
\draw [dashed] (4,4)--(2.5,4.5);
\fill (1,4) circle (2pt) node[above]{$v'$};
\fill (2.5,4.5) circle (2pt) node[above]{$v''$};
\fill (4,4) circle (2pt) node[above]{$v'''$};
\fill (-0.5,2) circle (2pt) node[below left]{$v_1$};
\draw (0.65,2) node {$\ldots$};
\fill (1.8,2) circle (2pt) node[below]{$v_{\frac{\Delta}{2}-1}$};
\fill (4.6,2) circle (2pt) node[below]{$v_{\frac{\Delta}{2}+2}$};
\fill (1.9,3.8) circle (2pt) node[below]{$v_{\frac{\Delta}{2}}$};
\fill (3.1,3.8) circle (2pt) node[below]{$v_{\frac{\Delta}{2}+1}$};
\draw (5.05,2) node {$\ldots$};
\fill (5.5,2) circle (2pt) node[below right]{$v_{\Delta}$};
\fill (-0.7,0) circle (2pt);
\fill (-0.3,0) circle (2pt);
\fill (5.3,0) circle (2pt);
\fill (5.5,0) circle (2pt);
\fill (5.7,0) circle (2pt);
\draw (1.9,3.8)--(2.5,4.5);
\draw (3.1,3.8)--(2.5,4.5);
\draw (1,4)--(-0.5,2);
\draw (1,4)--(1.8,2);
\draw (4,4)--(4.6,2);
\draw (4,4)--(5.5,2);
\draw (-0.5,2)--(-0.7,0);
\draw (-0.5,2)--(-0.3,0);
\draw (5.5,2)--(5.3,0);
\draw (5.5,2)--(5.5,0);
\draw (5.5,2)--(5.7,0);
\end{tikzpicture}
}
\end{center}
\caption{Construction from the proof of Theorem~\ref{thm:parallel} for the case $\Delta \equiv 0 \pmod{4}$}
\label{fig:construction4}
\end{figure}

Analogously as in the first case, note that in $G'$ all except the newly added vertices are of the same degree as in $G$, while $d_{G'}(v') = d_{G'}(v''') = \frac{\Delta}{2}$ and $d_{G'}(v'') =  4$ (the last two statements are true for all new vertices). It follows that $\Delta(G') < \Delta$. Since $\Delta \ge 6$, we then also infer that $\delta(G')\ge 4$. Because $\Delta \equiv 0 \pmod{4}$, the degrees $d_{G'}(v') = d_{G'}(v''') = \frac{\Delta}{2}$ are even, hence $G$ is Eulerian. By the induction assumption on $\Delta$, the graph $G'$ admits a parallel $2$-stable trace. Similarly as in previous case, if we use a path containing vertices $v'$, $v''$ and $v''$ as subtree $T$, it follows from repeated application Lemma~\ref{lem:contraction} that $G$ admits a parallel $2$-stable trace.

\vspace{5pt}

We have thus proved Theorem~\ref{thm:parallel}. 
\qed

\section{Alternative proof of Theorem~\ref{thm:parallel_d}}
\label{sec:parallel_d}

We now extend the results from previous section to present an alternative proof of Theorem~\ref{thm:parallel_d} (Theorem~$5.4$ from~\cite{fi-2013}).

\vspace{5pt}
\proof
Assume first that the graph $G$ admits a parallel $d$-stable trace. From Proposition~\ref{prop:stable} it follows that $\delta(G) > d$ for every graph $G$ that admits a $d$-stable trace. Assume that there exists an vertex $v$ of odd degree in $G$. Since every edge of a parallel double trace is used twice in the same direction, input and output degree of a parallel double trace $W$ would not match at $v$, which is absurd. Therefore it follows that $G$ is Eulerian and $\delta(G) > d$.

For the converse assume that graph $G$ is Eulerian and $\delta(G) > d$. Since $G$ is Eulerian, $\delta(G)$ is an even number. Furthermore, since for parallel $1$-stable traces and $2$-stable traces the theorem follows from Proposition~\ref{prp:parallel} and Theorem~\ref{thm:parallel}, respectively, we can assume that $d \geq 3$. Let $G'$ be a graph obtained from $G$ by replacing every vertex $v$ of degree $d_G(v) > 4$ with $(d_G(v)-2)/2$ new vertices, connected into a path $P_v$ and additionally connecting two endvertices of $P_v$ with three different neighbors of $v$ and each inner vertex of $P_v$ with two different remaining neighbors, so that each of the vertices from $N(v)$ is connected to exactly one vertex in $P_v$. It is not difficult to see that $G'$ is a $4$-regular graph and therefore by Theorem~\ref{thm:parallel} admits a parallel $2$-stable trace $W'$. Construct a parallel double trace $W$ in $G$ from $W'$ as follows. We start in an arbitrary vertex of $G'$ and follow $W'$. Let $a'=xy$ be an arc of $W'$ that we are currently traversing on our walk along $W'$. If for every $v$, $d_G(v) > 4$, $x,y \notin V(P_v)$, then we put $xy$ into $W$ so that the order of arcs from $W'$ is preserved. If for some $v$, $d_G(v) > 4$, $x \in P_v$ or $y \in P_v$, we replace $a'$ with $vy$ or $xv$, respectively. Finally, occurrences of the arcs with both endvertices contained in some $P_v$ are ignored in $W$.

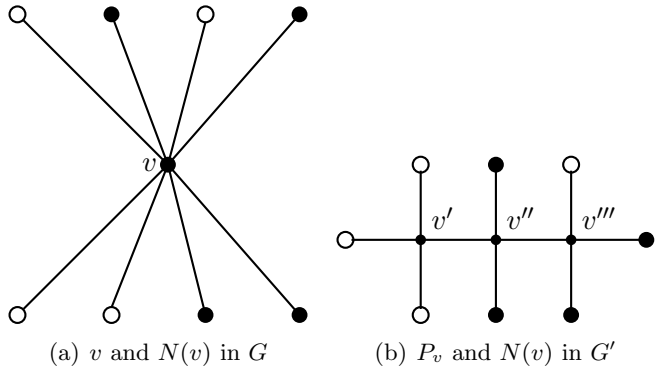
\begin{figure}[ht!]
\begin{center}
\subfigure[$v$ and $N(v)$ in $G$]
{
\begin{tikzpicture}[scale=1.0,style=thick]
\draw(2,2)--(0.08,0.08);
\draw(2,2)--(1.25,0.1);
\draw(2,2)--(2.5,0);
\draw(2,2)--(3.75,0);
\draw(2,2)--(0.08,3.92);
\draw(2,2)--(1.25,4);
\draw(2,2)--(2.5,3.9);
\draw(2,2)--(3.75,4);
\fill (2,2) circle (3pt) node[left]{$v$};
\draw (0,0) circle (3pt);
\draw (1.25,0) circle (3pt);
\fill (2.5,0) circle (3pt);
\fill (3.75,0) circle (3pt);
\draw (0,4) circle (3pt);
\fill (1.25,4) circle (3pt);
\draw (2.5,4) circle (3pt);
\fill (3.75,4) circle (3pt);
\end{tikzpicture}
}
\subfigure[$P_v$ and $N(v)$ in $G'$]
{
\begin{tikzpicture}[scale=1.0,style=thick]
\draw(1,1)--(0.08,1);
\draw(1,1)--(1,0.08);
\draw(1,1)--(1,1.92);
\draw(2,1)--(2,0);
\draw(2,1)--(2,2);
\draw(3,1)--(3,0);
\draw(3,1)--(3,1.92);
\draw(3,1)--(4,1);
\draw(1,1)--(2,1);
\draw(2,1)--(3,1);
\fill (1,1) circle (2pt) node[above right]{$v'$};
\fill (2,1) circle (2pt) node[above right]{$v''$};
\fill (3,1) circle (2pt) node[above right]{$v'''$};
\draw (0,1) circle (3pt);
\draw (1,0) circle (3pt);
\fill (2,0) circle (3pt);
\fill (3,0) circle (3pt);
\draw (1,2) circle (3pt);
\fill (2,2) circle (3pt);
\draw (3,2) circle (3pt);
\fill (4,1) circle (3pt);
\end{tikzpicture}
}
\end{center}
\caption{Structures of $N(v)$ in $G$ and $P_v$ in $G'$. Vertices contained in $X$ are colored black.}
\label{fig:coloring}
\end{figure}

We claim that the parallel double trace $W$ is a parallel $d$-stable trace of the graph $G$. We assume conversely and denote an arbitrary vertex in which $W$ has a repetition of order $\leq d$ with $v$. Denote the maximal order of $(\leq d)$-repetition at $v$ with $d'$. Since we used the same construction as in the proof of Theorem~\ref{thm:parallel}, it follows that $W$ is a parallel $2$-stable trace (and $d' > 2$). From the symmetry of repetitions it then also follows that $d_G(v) > d' + 2$, since otherwise $W$ would have at least one $1$-repetition or one $2$-repetition at $v$ (therefore also $d_G(v) \geq 8$). It is then also not difficult to see that every repetition in a parallel double trace is of even order. Let $X$ be a subset of $N(v)$ containing vertices from a maximal repetition at $v$ (note that $|X| = d'$). There exists a path $P_v$ in $G'$ that during the construction replaced $v$ from $G$. To make the argument more transparent, we imagine vertices from $P_v$ arranged in a horizontal line with all the neighbors of $v$ except for two, lying directly above or below vertices of $P_v$. The remaining two neighbors of vertex $v$ are aligned at the beginning and at the end of the horizontal line containing vertices from $P_v$. Fig.~\ref{fig:coloring} (b) shows $P_v$ with vertices from $N(v)$ in $G'$ for $d_G(v) = 8$ ($v'$, $v''$, and $v'''$ are the vertices replacing $v$ in $G'$). Next, we color vertices from $N(v)$ with two colors---black and white, so that vertices from $X$ are colored black while vertices from $N(v) \setminus X$ are colored white. Example of such a coloring can be seen at Fig.~\ref{fig:coloring}.

\begin{figure}[ht!]
\begin{center}
\subfigure[$b$ is the only black neighbor of $v'$]
{
\begin{tikzpicture}[scale=1.0,style=thick]
\draw (0.1,1)--(1,1);
\draw (1,0.1)--(1,1);
\draw (1,1)--(1,1.9);
\draw (1,1)--(2,1);
\draw (2,1)--(2,0.1);
\draw (2,1)--(2,2);
\draw (2,1)--(3,1);
\draw (3,1)--(3,2);
\draw (3,1)--(3,0);
\draw (3,1)--(4,1);
\draw (0,1) circle (3pt);
\fill (1,1) circle (2pt);
\draw (1,0) circle (3pt);
\draw (1,2) circle (3pt);
\fill (2,1) circle (2pt) node[above right]{$v'$};
\draw (2,0) circle (3pt);
\fill (2,2) circle (3pt) node[above right]{$b$};
\fill (3,1) circle (2pt) node[above right]{$v''$};
\fill[gray] (3,0) circle (2.9pt);
\fill[gray] (3,2) circle (2.9pt);
\draw (3,0) circle (3pt);
\draw (3,2) circle (3pt);
\fill[gray] (4,1) circle (2.9pt);
\draw (4,1) circle (3pt);
\end{tikzpicture}
}
\subfigure[Both neighbors of $v'$ from $N(v)$ are black]
{
\begin{tikzpicture}[scale=1.0,style=thick]
\draw (0.1,1)--(1,1);
\draw (1,0.1)--(1,1);
\draw (1,1)--(1,1.9);
\draw (1,1)--(2,1);
\draw (2,1)--(2,0);
\draw (2,1)--(2,2);
\draw (2,1)--(3,1);
\draw (3,1)--(3,2);
\draw (3,1)--(3,0);
\draw (3,1)--(4,1);
\draw (0,1) circle (3pt) node[above right]{$w_1$};
\fill (1,1) circle (2pt);
\draw (1,0) circle (3pt) node[above right]{$w_3$};
\draw (1,2) circle (3pt) node[above right]{$w_2$};
\fill (2,1) circle (2pt) node[above right]{$v'$};
\fill (2,0) circle (3pt) node[above right]{$b'$};
\fill (2,2) circle (3pt) node[above right]{$b$};
\fill (3,1) circle (2pt) node[above right]{$v''$};
\fill[gray] (3,0) circle (2.9pt);
\fill[gray] (3,2) circle (2.9pt);
\draw (3,0) circle (3pt);
\draw (3,2) circle (3pt);
\fill[gray] (4,1) circle (2.9pt);
\draw (4,1) circle (3pt);
\end{tikzpicture}
}
\end{center}
\caption{Two cases of the structure of $P_v$ (of $v'$ and $b$ to be more precise). Vertices for which the color is not determined are colored grey.}
\label{fig:cases1}
\end{figure}
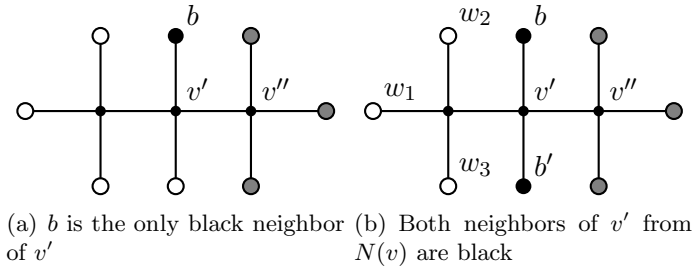

Since the subset $N(v) \setminus X$ is also a repetition, the arguments used hereinafter are true for black and white vertices and we can, without loss of generality, assume that the neighbor of $N(v)$, lying farmost to the left in the above mentioned horizontal line is colored white. We next move along this horizontal line and denote the first black vertex that we meet (below or above the line) with $b$. Denote its neighbor in $P_v$ with $v'$. Since there are at least four black vertices, $v'$ is not the farmost right vertex from $P_v$. Therefore, we can also denote the right neighbor of $v'$ from $P_v$ with $v''$ and consider two cases. In the first case $b$ is the only neighbor of $v'$ ($\notin P_v$) colored black (Fig.~\ref{fig:cases1} (a)), while in the second case also the second neighbor of $v'$ ($\notin P_v$) is colored black (Fig.~\ref{fig:cases1} (b)). In both cases we can, without loss of generality, assume that an edge $bv'$ is traversed twice in the direction toward $v'$ in $W'$ (that is, arc $bv'$ is traversed twice in $W$, while arc $v'b$ does not appear in $W'$). The fact that $W$ has an $X$-repetition implies that every time double trace $W$ comes to $v$ from a vertex in $X$ it also continues to a vertex in $X$ and, consequently, that every time a double trace $W'$ comes to a vertex in $P_v$ from a (black colored) vertex in $X$ it also leaves to a (black colored) vertex in $X$. Note that in between $W'$ can traverse other vertices from $P_v$ and that this applies for all appearances of verb {\em continue} from now on until the end of this section. Analogously is true for (white colored) vertices from $N(v) \setminus X$. Therefore, in $W'$ there exist two subsequences which start with $bv'$, continue on some other vertices from $P_v$ and end in two from $b$ different vertices from $X$. In the first case, when $b$ is the only black colored neighbor of $v'$, the subsequence $b \rightarrow v' \rightarrow v''$ has to appear twice in $W'$, since otherwise $W'$ can not continue (twice) from $b$ to a black colored vertex without previously traversing white vertex. This contradicts the fact that $W'$ is a parallel $2$-stable trace, since $bv'v''$ is a $2$-repetition at $v'$. In the second case, we denote the set of white vertices that appear to the left of $b$ with ${\cal W} = \{w_1, \ldots, w_l\}$. (Note that $l$ is an odd integer.) For an example, see Fig.~\ref{fig:cases1} (b), where those vertices are denoted with $w_1, w_2$, and $w_3$. Next, we denote the second black colored neighbor of $v'$ from $N(v)$ with $b'$. The subsequence $b \rightarrow v' \rightarrow b'$ can appear at most once in $W'$ (otherwise $W'$ would have a $2$-repetition at $v'$). Assume next that for every $w \in {\cal W}$, $w$ continues to a vertex in ${\cal W}$. Then vertices from ${\cal W}$ form an odd repetition in $W'$, which can not appear in a parallel $2$-stable trace. Therefore, at least one vertex $w$ from ${\cal W}$ has to continue to a white colored vertex not included in ${\cal W}$ (that is, $w$ continues to a white colored vertex to the right of $b$). If subsequence $b \rightarrow v' \rightarrow b'$ does not appear in $W'$ it follows that edge $v'v''$ (arc $v'v''$ and $v''v'$) is used more than twice in $W'$: at least once to connect a vertex from ${\cal W}$ to a white colored vertex not in ${\cal W}$, twice to connect $b$ to a (black colored) vertex in $X$ different from $b'$, and twice to connect $b'$ to a (black colored) vertex in $X$ different from $b$, which is absurd. If subsequence $b \rightarrow v' \rightarrow b'$ does appear in $W'$ it (in addition to multiple appearances of $v'v''$) follows that edge $v'v''$ is not parallel in $W'$. Since all the black colored vertices except $b$ and $b'$ are to the right of $v'$ both $b \rightarrow v' \rightarrow v''$ and $v'' \rightarrow v' \rightarrow b'$ have to appear in $W'$, which is also absurd.

Since $v$ was an arbitrary vertex in $G$ and $d'$ was an arbitrary integer, $2 < d' \leq d$, it follows that $W$ is a parallel $d$-stable trace of $G$ and therefore Theorem~\ref{thm:parallel_d} is proved.
\qed

\section{Concluding remarks}
\label{sec:conclusion}

In this section we present two concepts which we assumed could be used for constructing parallel $2$-stable traces. Unfortunately, it has turned out, when proving Theorem~\ref{thm:parallel}, that there exist graphs admitting only parallel $2$-stable traces which can not be realized using the here described constructions.

The first construction goes as follows. Let $G$ be an Eulerian graph with $n$ vertices (denoted with $v_1, \ldots, v_n$) fulfilling conditions of Theorem~\ref{thm:parallel} and let $W'$ be an Eulerian circuit of $G$. $W'$ induces a set of functions $\Pi' = \{\pi'_1, \ldots, \pi'_n\}$, where $\pi_i': N(v_i) \longrightarrow N(v_i)$, $\pi'_i(v) = u$ if and only if $v \rightarrow v_i \rightarrow u$ or $u \rightarrow v_i \rightarrow v$ are sequences in $W'$, for $1 \leq i \leq n$. Note that $u \neq v$, because $G$ is simple and $W'$ traverses every edge exactly once. Suppose that $W''$ is another Eulerian circuit in $G$ such that $W''$ induces a set of functions $\Pi'' = \{\pi''_1, \ldots, \pi''_n\}$ with above described characteristics. In addition demand that edges are traversed in the same direction as in $W'$, and that if $\pi'_i(v) = u$ then $\pi''_i(v) \neq u$ and $\pi''_i(u) \neq v$. Concatenate Eulerian circuits $W'$ and $W''$ into a double trace $W$ in an arbitrary vertex $v$. It is obvious from the construction that every edge is traversed twice in the same direction in $W$ and that $W$ is without $1$-repetitions and $2$-repetitions in any vertex other than $v$. Hence, if a graph $G$ admits two Eulerian circuits with above described characteristics, then $G$ admits parallel $2$-stable trace as well. 

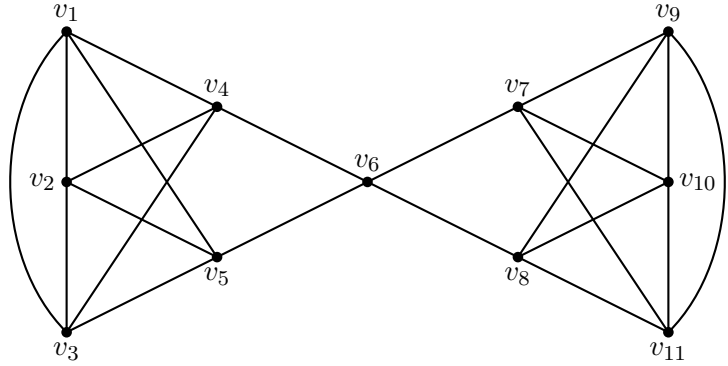
\begin{figure}[ht!]
\begin{center}
\begin{tikzpicture}[scale=1,style=thick]
\fill (1,0) circle (2pt) node[below]{$v_3$};
\fill (1,2) circle (2pt) node[left]{$v_2$};
\fill (1,4) circle (2pt) node[above]{$v_1$};
\fill (3,1) circle (2pt) node[below]{$v_5$};
\fill (3,3) circle (2pt) node[above]{$v_4$};
\fill (5,2) circle (2pt) node[above]{$v_6$};
\fill (7,1) circle (2pt) node[below]{$v_8$};
\fill (7,3) circle (2pt) node[above]{$v_7$};
\fill (9,0) circle (2pt) node[below]{$v_{11}$};
\fill (9,2) circle (2pt) node[right]{$v_{10}$};
\fill (9,4) circle (2pt) node[above]{$v_9$};
\draw (5,2)--(3,1);
\draw (5,2)--(3,3);
\draw (3,1)--(1,0);
\draw (3,1)--(1,2);
\draw (3,1)--(1,4);
\draw (3,3)--(1,0);
\draw (3,3)--(1,2);
\draw (3,3)--(1,4);
\draw (1,0)--(1,2);
\draw (1,2)--(1,4);
\draw (1,0)  .. controls (0,1) and (0,3) .. (1,4);
\draw (5,2)--(7,1);
\draw (5,2)--(7,3);
\draw (7,1)--(9,0);
\draw (7,1)--(9,2);
\draw (7,1)--(9,4);
\draw (7,3)--(9,0);
\draw (7,3)--(9,2);
\draw (7,3)--(9,4);
\draw (9,0)--(9,2);
\draw (9,2)--(9,4);
\draw (9,0)  .. controls (10,1) and (10,3) .. (9,4);
\end{tikzpicture}
\end{center}
\caption{Graph whose parallel $2$-stable traces cannot be constructed by concatenating two Eulerian circuits}
\label{fig:concatenating}
\end{figure}

It turns out that we cannot always construct a parallel $2$-stable trace of $G$ by concatenating two Eulerian circuits. For instance, the graph $G$ from Fig.~\ref{fig:concatenating} has a parallel $2$-stable trace: $v_1 \rightarrow v_2 \rightarrow v_3  \rightarrow v_1  \rightarrow v_2  \rightarrow v_4  \rightarrow v_1 \rightarrow v_5 \rightarrow v_2 \rightarrow v_3 \rightarrow v_4 \rightarrow v_6 \rightarrow v_5 \rightarrow v_2  \rightarrow v_4  \rightarrow v_6 \rightarrow v_7 \rightarrow v_9 \rightarrow v_8 \rightarrow v_6 \rightarrow v_7 \rightarrow v_{10} \rightarrow v_8 \rightarrow v_{11} \rightarrow v_7 \rightarrow v_9 \rightarrow v_{10} \rightarrow v_{11} \rightarrow v_7 \rightarrow v_{10} \rightarrow v_{11} \rightarrow v_9 \rightarrow v_8 \rightarrow v_{11} \rightarrow v_9 \rightarrow v_{10} \rightarrow v_8 \rightarrow v_6 \rightarrow v_5 \rightarrow v_3 \rightarrow v_1 \rightarrow v_5 \rightarrow v_3 \rightarrow v_4 \rightarrow v_1$, but because of the cut vertex $v_6$, from any Eulerian circuit $W$ of $G$ we cannot construct another Eulerian circuit with the described properties.     

The main idea of the second construction is to find a parallel $2$-stable trace in each block of a graph $G$ and then concatenate them into a parallel $2$-stable trace of the graph $G$. Let again $G$ be an Eulerian graph fulfilling the conditions of Theorem~\ref{thm:parallel}. Denote blocks of $G$ with $B_1, \ldots, B_k$ and cutvertices with $v_1, \ldots, v_l$. Find first a parallel $2$-stable trace $W_i$ in block $B_i$. Concatenate parallel $2$-stable traces into a parallel $2$-stable trace of $G$ in corresponding cutvertices. When concatenating, one has to be careful that no $1$-repetitions and $2$-repetitions appear.

Similar as for the first construction, none of the parallel $2$-stable traces of the graph $G$ from Fig.~\ref{fig:concatenating} can not be constructed by concatenating parallel $2$-stable traces in its blocks. Vertex $v_6$ is a unique cutvertex of the graph $G$ and it is contained in both of its blocks. Since $v_6$ is of degree $2$ in both blocks of the graph $G$, none of them admit parallel $2$-stable trace.  Similar problem occurs if one or more blocks of $G$ are bridges. 

Next possible improvement could instead of parallel $2$-stable traces in blocks demand parallel $1$-stable traces where $2$-repetitions (or $1$-repetitions if block is a bridge) would be allowed at cutvertices but are then later removed during the concatenation into a parallel $2$-stable trace of the whole graph.

An attempt to find efficient algorithms for constructing and counting stable traces of graphs was made in~\cite{bas-2017}. It would be of interest to characterize graphs that do not have any of the two above described properties of graphs from Fig.~\ref{fig:concatenating} 
and then try to improve the algorithms from~\cite{bas-2017} by using the above described constructions for those special cases of graphs.

\section*{Acknowledgements}

The author is grateful to Sandi Klav\v zar and anonymous reviewers for several significant remarks and suggestions which were of great help. The authors acknowledge the financial support from the Slovenian Research Agency H2020 SME2 and the SPIRIT Slovenia - Public Agency for Entrepreneurship, Internationalization, Fore­ign Investments and Technology - KKIPP.


\end{document}